\documentclass{amsart} 
\numberwithin{equation}{section}
\newtheorem{theorem}[equation]{Theorem}
\newtheorem*{theorem*}{Theorem}
\newtheorem{proposition}[equation]{Proposition}
\newtheorem{corollary}[equation]{Corollary}

\newtheorem{lemma}[equation]{Lemma}

\theoremstyle{remark}
\newtheorem{remark}[equation]{Remark}

\newtheorem*{remark*}{Remark}
\newtheorem*{note}{Note}

\newtheorem{example}[equation]{Example}
\newtheorem*{example*}{Example}

\DeclareMathOperator{\rk}{rk}
\DeclareMathOperator{\ad}{ad}
\DeclareMathOperator{\End}{End}

\begin{document}
\title[Stable bundles on Hopf manifolds]
{Stable bundles on Hopf manifolds}

\author{Ruxandra Moraru}
\address{Department of Mathematics, University of Toronto, 100 St George Street, 
Toronto, Ontario, Canada, M5S 3G3}
\email{moraru@math.toronto.edu}


\begin{abstract}
In this paper, we study holomorphic vector bundles on 
(diagonal) Hopf manifolds. 
In particular, we give a description of moduli spaces of stable bundles
on generic (non-elliptic) Hopf surfaces.
We also give a classification of stable
rank-2 vector bundles on generic Hopf manifolds of complex dimension 
greater than two. 
\end{abstract}

\maketitle

\section{Introduction}
In this paper, we study the stability properties of
holomorphic vector bundles on diagonal Hopf manifolds.
Recall that a Hopf manifold is defined as the quotient of the punctured $n$-space 
$\mathbb{C}^n \backslash \{ 0 \}$ by an infinite cyclic group, generated 
by a contraction of $(\mathbb{C}^n,0)$. If the contraction
is multiplication by a diagonal matrix, then the Hopf manifold
is called {\em diagonal}. All Hopf manifolds are non-algebraic.
In particular, every diagonal Hopf manifold is diffeomorphic to
$S^1 \times S^{2n+1}$, implying that it is non-K\"{a}hlerian. 
A generic Hopf manifold possesses very few curves, the only ones being $n$
elliptic curves corresponding to the coordinate axes in $\mathbb{C}^n$.
But there exist, nevertheless, Hopf manifolds with infinite families of curves. For example, if the
contraction defining the manifold is a multiple of the identity, then the manifold admits
an elliptic fibration. 

Holomorphic vector bundles on elliptically fibred Hopf manifolds are by now well-understood.
In the case of surfaces, these bundles have been completely classified and
a detailed analysis of their moduli spaces can be found in \cite{Moraru,Brinzanescu-Moraru2}.
Moreover, for Hopf manifolds of dimension greater than two, the question
has been settled by Verbitsky in \cite{V1}, where he studies
stable bundles on positive principal elliptic fibrations.
More generally, he proves that on a (possibly non-elliptic) 
diagonal Hopf manifold of dimension greater than two,
any coherent sheave $\mathcal{F}$ is {\em filtrable}, that is, 
admits a filtration by a sequence of coherent sheaves
\[ 0 = \mathcal{F}_0 \subset \mathcal{F}_1 \subset \cdots \subset 
\mathcal{F}_r = \mathcal{F} \]
with $\rk \mathcal{F}_i / \mathcal{F}_{i-1} \leq 1$ \cite{V1,V2}.

Less is known about the classification and stability properties of
bundles on generic (non-elliptic) Hopf manifolds. 
Although some partial results have been obtained in this direction by Mall 
\cite{Mall-filtrable,Mall-rank-2}, his study has focused only on vector bundles 
whose pullback to the universal cover $\mathbb{C}^n \backslash \{ 0 \}$ is 
holomorphically trivial (such bundles are given by factors of automorphy).  
In particular, he proves that the pullback of a vector bundle on a Hopf manifold 
to $\mathbb{C}^n \backslash \{ 0 \}$ is holomorphically trivial if and only if it 
possesses a filtration by vector bundles, giving in the process a partial classification 
such bundles on generic Hopf manifolds. 
However, as shown in this paper, most vector bundles on diagonal
Hopf manifolds do not admit filtrations by vectors bundles; in fact, on Hopf
surfaces, vector bundles are generically non-filtrable. 

The paper is organised as follows. 
We begin by recalling some definitions and topological properties of Hopf manifolds;
holomorphic vector bundles are described  in sections three and four.  
For simplicity, we restrict our presentation on generic Hopf manifolds
to rank-2 vector bundles; nevertheless, similar results holds for vector bundles of arbitrary
rank. 
The third section of the paper is devoted to the study of bundles on surfaces.
We begin by proving that, on a diagonal Hopf surface, holomorphic
vector bundles possess a filtration by vector bundles
if and only if they are topologically trivial; bundles with non-trivial 
second Chern classes are, however, generically non-filtrable and therefore stable.
We then give a classification of stable filtrable bundles.
Moduli spaces of stable bundles on diagonal Hopf surfaces 
admit natural Poisson structures; we describe their associated symplectic leaves 
for elliptic Hopf surfaces. Note that in the elliptic case, 
the moduli also admit the structure of algebraically
completely integrable Hamiltonian systems \cite{Moraru}.
Finally, stability conditions for vector bundles on generic Hopf
manifolds of dimension greater than two are given in the last section; 
in particular, we show that there exist stable bundles on them.
\vspace{.1in}

{\bf Acknowledgements.} 
The author would like to express her gratitude to Misha Verbitsky 
for explaining to her the filtrability of vector bundles
on higher dimensional Hopf manifolds
as well as for valuable discussions and suggestions. 
She would also like to thank 
Jacques Hurtubise, Boris Khesin, and Vasile Br\^{\i}nz\u{a}nescu
for useful comments, and the
Department of Mathematics at the University of Glasgow 
for their hospitality during the preparation of part of this article.

\section{Preliminaries}
\label{preliminaries}

A diagonal Hopf manifold $X$
is defined as the quotient of $\mathbb{C}^n_\ast := \mathbb{C}^n \backslash \{ 0 \}$ 
by the cyclic group generated by a contraction of $(\mathbb{C}^n,0)$
of the form
\[ \begin{array}{rcl}
\mu: \mathbb{C}^n & \longrightarrow & \mathbb{C}^n \\ 
(z_1, \dots, z_n) & \mapsto & (\mu_1 z_1, \dots, \mu_n z_n),
\end{array} \]
where $\mu_1, \dots, \mu_n$ are complex numbers such that 
$0 < | \mu_1 | \leq | \mu_2 | \leq \dots \leq | \mu_n | < 1$.
Every diagonal Hopf manifold is diffeomorphic to $S^1 \times S^{2n+1}$.
Moreover, its Hodge numbers are all zero except 
for $h^{0,0}_X = h^{0,1}_X =  h^{n+1,n}_X = h^{n+1,n+1}_X = 1$ \cite{Mall-line bundles}.

\subsection{Notation}

We begin by fixing some notation.

\begin{itemize}
\item
Denote by $p: \mathbb{C}^n_\ast \rightarrow X$ the canonical projection map.

\item
$U_i : = \{ z_i \neq 0 \}$ is an open subset of $\mathbb{C}^n_\ast$ for all $i = 1, \dots, n$.

\item
$X_i := p(U_i)$ is an open subset of $X$ for all $i = 1, \dots, n$.

\item
$H_i := p(\{ z_i = 0 \})$ is a hypersurface in $X$, for all $i = 1, \dots, n$, that is isomorphic to
the Hopf manifold of dimension $n-1$ corresponding to the diagonal matrix
$(\mu_1, \dots, \hat{\mu_i}, \dots, \mu_n)$.

\item
$H_{k_i k_j} := p(\{ z_i^{k_i} = z_j^{k_j} = 0 \})$, for
$k_i,k_j \geq 0$ and $0 \leq i,j \leq n$ such that $i \neq j$, 
is a codimension 2 subvariety of $X$. 

\item
$T_i := p(\{ (0, \dots, z_i, \dots, 0) \})$
is the elliptic curve $\mathbb{C}^\ast / \mu_i$ for all $i = 1, \dots, n$.
Note that $H_i \cap T_i = \emptyset$ for all $i = 1, \dots, n$. 
Moreover, $T_i \cap T_j = \emptyset$, if $i \neq j$.
\end{itemize}

\subsection{Classical Hopf manifolds}

A diagonal Hopf manifold is called {\em classical} if $\mu_1 = \mu_2 = \cdots = \mu_n = \mu$.
These manifolds admit a natural holomorphic elliptic fibration 
\[ \begin{array}{rcl}
\pi: X & \rightarrow & \mathbb{P}^n \\
(z_1, \cdots, z_n) & \mapsto & [ z_1 : \cdots : z_n],
\end{array} \]
with fibre the elliptic curve $T = \mathbb{C}^\ast/\mu$.
Moreover, the relative Jacobian of $X \stackrel{\pi}{\rightarrow} \mathbb{P}^n$ 
is isomorphic to
\[ J(X) = \mathbb{P}^n \times T^\ast \stackrel{p_1}{\rightarrow} \mathbb{P}^n, \]
where $T^\ast$ denotes the dual elliptic curve determined by a non-canonical identification 
$T^\ast := \text{Pic}^0(T) \cong T$. 

\subsection{Generic Hopf manifolds}

A diagonal Hopf manifold is called {\em generic} if there are no non-trivial relations between 
the $\mu_i$'s of the form 
\[ \prod_{i \in A} \mu_i^{r_i} = \prod_{j \in B} \mu_j^{r_j}, \]
where $r_i,r_j \in \mathbb{N}$, $A \cap B = \emptyset$, and $A \cup B = \{ 1, \dots, n \}$.
It is important to note that a generic Hopf manifold only contains $n$ irreducible curves, namely
the images $T_1, \dots, T_n$ of the punctured $z_1$-, $\dots$, $z_n$-axes 
(see \cite{BPV} for the case $n=2$). 
Although these curves are elliptic, the manifold itself does not
admit an elliptic fibration. 
Moreover, given that there are no relations between the $\mu_i$'s, the $H_i$'s are the 
only irreducible hypersurfaces in $X$.

\subsection{Hopf surfaces}
\label{non-primary-Hopf}
Diagonal Hopf surfaces can be divided into four categories: classical, generic, resonant,
and hyperresonant. The last two are defined by diagonals $(\mu_1, \mu_2)$
such that $\mu_1^p = \mu_2^q$ for some integers $p$ and $q$
(the resonant case corresponds to $p=1$).
Note that a (hyper)resonant surface $X$ admits an elliptic fibration 
$\pi: X \rightarrow \mathbb{P}^1$ with singular fibres 
the curve(s) $T_1$ (and $T_2$). 
However, it can be covered by the 
classical Hopf surface given by diagonal $(\mu_2^{1/p},\mu_2^{1/p})$,
whose elliptic fibration does not have singular fibres.
Vector bundles on the (hyper)resonant surface $X$ can therefore be studied by analysing
their pullback to this classical Hopf surface 
\cite{Brinzanescu-Moraru,Brinzanescu-Moraru1,Brinzanescu-Moraru2}.

In general, Hopf surfaces are defined as compact complex surfaces that admit
$\mathbb{C}^\ast$ as a universal covering.
Although every diagonal Hopf surface is diffeomorphic to $S^1 \times S^3$,
there exist many Hopf surfaces that are not.
Examples of such Hopf surfaces can be constructed  as follows.
For any integer $d$, let $\Theta_d^\ast$ denote the total space
of the line bundle $\mathcal{O}_{\mathbb{P}^1}(d)$ minus the zero section;
with respect to this notation, we have $\mathbb{C}^2_\ast \cong \Theta_{-1}^\ast$.
We cover $\Theta_d^\ast$ by the open sets $U_0$ and $U_{\infty}$
with coordinates $(z,t)$ and $(\xi,s) = (z^{-1},z^{-d} t)$, respectively.
Given complex number $\mu_1$ and $\mu_2$, with $0 < | \mu_1 | \leq | \mu_2 | < 1$,
and $\beta$ such that $\beta^d = \mu_1\mu_2^{-1}$, we define the following
$\mathbb{Z}$-action on $\Theta_d^\ast$: $(z,t) \mapsto (\beta z,\mu_1 t)$ on $U_0$,
and $(\xi,s) \mapsto (\beta^{-1}\xi,\mu_2 s)$ on $U_\infty$.
The quotient 
\[ X := \Theta_d^\ast / (\beta,\mu_1)\]
is a Hopf surface, which is non-primary when $d \neq -1$.
Note that if $d=-1$, then we simply have a diagonal Hopf surface:
\[ \mathbb{C}_\ast^2/(\mu_1,\mu_2) \cong \Theta_{-1}^\ast / (\mu_1^{-1}\mu_2,\mu_1).\]
Moreover, if $\mu_1 = \mu_2$ and $\beta = 1$, then $X$ is a principal elliptic 
fibration over $\mathbb{P}^1$ with fibre $\mathbb{C}^\ast/\mu_1$.
These Hopf surfaces are diffeomorphic to $S^1 \times C_d$, where
$C_d$ is the $S^1$-bundle over $S^2$ with Chern class $d$.
Consequently, we have $H^i(X,\mathbb{Z}) = \mathbb{Z}$, for $i = 0,1,4$, 
$H^2(X,\mathbb{Z}) = \mathbb{Z}_{|d|}$,
and $H^3(X,\mathbb{Z}) = \mathbb{Z} \oplus \mathbb{Z}_{|d|}$.

\subsection{Line bundles}
\label{line bundles}

The only divisors that exist on a generic Hopf manifold $X$
are linear combinations of the hypersurfaces $H_1, \dots, H_n$:
\[ \text{Div}(X) = \{ m_1 H_1 + \dots + m_n H_n \} = 
\mathbb{Z} \oplus \dots \oplus \mathbb{Z}, \]
with the relations $H_i \cdot H_j = 0$;
the canonical divisor is $K_X = -H_1 - \dots - H_n$.
Whereas divisors on a classical Hopf manifold $X$ are pullbacks of
hypersurfaces on $\mathbb{P}^{n-1}$.
In particular, the canonical divisor is $K_X = \pi^\ast K_{\mathbb{P}^{n-1}}$.

Although Hopf manifolds have few divisors,
there are many line bundles on them. For example, on any diagonal
Hopf manifold $X$, we have 
$\text{Pic}(X) = \text{Pic}^0(X) = \mathbb{C}^\ast$: line bundles correspond to
constant factors of automorphy. 
The line bundle given by the factor $a \in \mathbb{C}^\ast$ is constructed by
taking the quotient of the trivial line bundle $\bar{\mathbb{C}}$ on 
$\mathbb{C}^n_\ast$ by the following $\mathbb{Z}$-action:
\[ \begin{array}{rcl}
\mathbb{C}^n_\ast \times \mathbb{C} & \rightarrow & 
\mathbb{C}^n_\ast \times \mathbb{C}\\
(z, t) & \mapsto & (\mu z, at) .
\end{array} \]

\noindent
From now on, we shall denote by $L_a$ 
the line bundle corresponding to the factor $a$.
Note that the restriction of $\mathcal{O}_X(H_i)$ to the elliptic curve 
$T_i = \mathbb{C}^\ast / \mu_i$ is trivial, 
so that $\langle \mathcal{O}_X(H_i) \rangle = \mathbb{Z}$ is in the kernel of the 
natural restriction map $\text{Pic}(X) \stackrel{r}{\rightarrow} \text{Pic}^0(T_i) \rightarrow 0$.
Consequently, we have 
$\mathcal{O}_X(H_i) = L_{\mu_i}$ for all $i = 1, \dots, n$.

The cohomology of line bundles on classical and generic Hopf manifolds
the following \cite{Mall-line bundles}.
Given a line bundle $L_a$ on the Hopf manifold $X$, we denote
\[ h^{p,0} := h^p(X,L_a). \]
For a classical Hopf manifold $X$, 
given by the diagonal $(\mu, \dots, \mu)$, we have:
\[ h^{p,0} = \left\{ \begin{array}{cl}
h^p(\mathbb{P}^n,\mathcal{O}_{\mathbb{P}^n}(m)) & \mbox{if $a = \mu^m$ for some integer $m$},\\
0 & \mbox{otherwise}.
\end{array} \right. \]

We now consider a generic Hopf manifold $X$, 
given by the diagonal $(\mu_1, \dots, \mu_n)$.
If $n=2$, then given the line bundle $L_a$ on the generic Hopf surface $X$, we have:
\[ \begin{array}{l}
h^{0,0} = \left\{ \begin{array}{cl}
1 & \mbox{if $a = \mu_1^{m_1} \mu_2^{m_2}$, with $m_1, m_2 \geq 0$},\\
0 & \mbox{otherwise};
\end{array} \right. \\\\

h^{1,0} = h^{0,0} + h^{2,0}; \\\\

h^{2,0} = \left\{ \begin{array}{cl}
1 & \mbox{if $a = \mu_1^{m_1} \mu_2^{m_2}$, with $m_1, m_2 < 0$},\\
0 & \mbox{otherwise}.
\end{array} \right. 
\end{array}\]

\noindent
Finally, for $n \geq 3$, the cohomology groups of the line bundle 
$L_a$ are the following:  
\[ \begin{array}{l}
h^{0,0} = h^{1,0} = \left\{ \begin{array}{cl}
1 & \mbox{if $a = \mu_1^{m_1} \dots \mu_n^{m_n}$, with $m_1, \dots, m_n \geq 0$},\\
0 & \mbox{otherwise};
\end{array} \right. \\\\

\mbox{$h^{p,0} = 0$ if $2 \leq q \leq n-2$}; \\\\

h^{n-1,0} = h^{n,0} = \left\{ \begin{array}{cl}
1 & \mbox{if $a = \mu_1^{m_1} \dots \mu_n^{m_n}$, with $m_1, \dots, m_n < 0$},\\
0 & \mbox{otherwise}.
\end{array} \right. 
\end{array} \]

\subsection{Finite coverings}

In this section, we describe (finite) cyclic coverings of Hopf surfaces.
We begin by noting that elliptically fibred Hopf surfaces admit many cyclic coverings
as they correspond to pullbacks of cyclic coverings on $\mathbb{P}^1$ via the 
projection $\pi$. This is certainly not the case for generic Hopf surfaces. 
Let us consider a generic Hopf surface $X$ given by the diagonal $(\mu_1,\mu_2)$.
Let $\varphi: Y \rightarrow X$ be the $r$-cyclic covering of $X$ branched along the
smooth divisor $B$ on $X$ and determined by the line bundle $\mathcal{L} \rightarrow X$, 
where $\mathcal{O}_X(B) = \mathcal{L}^{\otimes r}$ (see \cite{BPV,F}).
The smooth variety $Y$ is then defined as the subvariety of the total space of the line 
bundle $\mathcal{L}$ given by
\[ Y := \{ y  \ : \ y^r = s \},\]
where $s$ is a global section of $\mathcal{L}^{\otimes r}$ with effective divisor $(s) = B$.
However, since $X$ is a generic Hopf surface,
we have $h^0(X,\mathcal{O}_X(B)) = 1$, 
implying that there is a unique cyclic covering corresponding 
to each $\mathcal{L}$.

\begin{lemma}\label{double covers}
Suppose that $X = \mathbb{C}^2_\ast/(\mu_1,\mu_2)$ is a generic Hopf surface. 
An $r$-cyclic covering 
of $X$ is then either a generic Hopf surface or a non-elliptic non-primary Hopf surface
of the form $\Theta^\ast_{-r}/ (\beta,\mu_1)$, with $\beta^r = \mu_1$ 
(see section \ref{non-primary-Hopf} for notation). 
\end{lemma}
\begin{proof}
On a generic Hopf surface, there are only four possibilities for the effective reduced 
divisor $B$, namely, $B = 0, \ T_1, \ T_2$, or $T_1 + T_2$, with associated line bundles
$\mathcal{O}_X(B) = \mathcal{O}_X, \ L_{\mu_2}, \ L_{\mu_1}$, and $L_{\mu_1\mu_2}$,
respectively. We treat each case seperately.

If $B = 0$, then the line bundle $\mathcal{L}$ is an $r$-th root of unity of order $k$, 
where $k$ is an integer that divides $r$. 
If $k=0$, that is, $\mathcal{L} = \mathcal{O}_X$, then $Y$ is a disconnected surface made up
of $r$ copies of $X$. If $\mathcal{L}$ has instead order $r$, then $Y$ is a connected smooth
surface that is an unramified $r$-to-one covering of $X$. In particular, there is a unique 
unramified double cover that is the Hopf
surface $Y := \mathbb{C}^2_\ast / (\mu_1^2,\mu_2^2)$ with projection onto $X$ given
by $(z_1,z_2) \mapsto (z_1,z_2)$.
Finally, if the order of $\mathcal{L}$ is $k \neq 0,m$, then $m = kl$ for some integer $l$
and $Y$ is disconnected surface made up of $l$ copies of an unramified $k$-to-one cover 
$Y'$ of $X$.

If $B = T_1$, then the line bundle $\mathcal{L}$ is
given by a factor of automorphy $\alpha \in \mathbb{C}^\ast$ such that $\alpha^r = \mu_2$.
The induced $r$-to-one cover of $X$ is then the Hopf
surface $Y := \mathbb{C}^2_\ast / (\mu_1,\alpha)$ with projection onto $X$
given by $(z_1,z_2) \mapsto (z_1,z_2^r)$. Similarly,
one sees that the $r$-to-one cover of $X$ determined by $B = T_2$ are Hopf
surfaces of the form $Y := \mathbb{C}^2_\ast / (\alpha,\mu_2)$ with $\alpha^r = \mu_1$.

Finally, consider $B = T_1 + T_2$; 
choose $\beta \in \mathbb{C}^\ast$ such that $\beta^r = \mu_1^{-1}\mu_2$.
Using the notation of section \ref{non-primary-Hopf}, 
we have $X = \Theta_{-1}^\ast / (\mu_1^{-1}\mu_2,\mu_1)$ and its 
$r$-to-one cover is the Hopf
surface $Y :=  \Theta_{-m}^\ast / (\beta, \mu_1)$ with projection onto $X$
given by $(z,t) \mapsto (z^r,t)$. 
\end{proof}

\subsection{Degree and stability}
\label{degree}
\indent
The degree of a vector bundle can be defined on any compact
complex manifold $M$. Let $d = \dim_\mathbb{C} M$. 
A theorem of Gauduchon's \cite{Gauduchon} states that any hermitian metric on $M$ 
is conformally equivalent to a metric, 
called a {\it Gauduchon metric}, whose associated (1,1) form 
$\omega$ satisfies $\partial \bar{\partial} \omega^{d-1} = 0$.
Suppose that $M$ is endowed with such a metric and let $L$ be a holomorphic 
line bundle on $M$. The {\it degree of $L$ with respect to $\omega$} is 
defined \cite{Buchdahl}, up to a constant factor, by
\[ \deg{L} := \int_M F \wedge \omega^{d-1},\]
where $F$ is the curvature of a hermitian connection on $L$, compatible with 
$\bar{\partial}_L$. Any two such forms $F$ differ by a 
$\partial \bar{\partial}$-exact form. Since 
$\partial \bar{\partial} \omega^{d-1} = 0$, 
the degree is independent of the choice of connection and is therefore
well defined. This notion of degree is an extension of the K\"{a}hler case. 
If $M$ is K\"{a}hler, we get the usual topological degree defined
on K\"{a}hler manifolds; but in general, this degree is 
not a topological invariant, for it 
can take values in a continuum (see below).

Having defined the degree of holomorphic line bundles, we define the 
{\it degree} of a torsion-free coherent sheaf $\mathcal{E}$ on $M$ by
\[ \deg(\mathcal{E}) := \deg(\det{\mathcal{E}}), \]
where $\det{\mathcal{E}}$ is the determinant line bundle of $\mathcal{E}$,
and the {\it slope of $\mathcal{E}$} by 
\[ \mu(\mathcal{E}) := \deg(\mathcal{E})/\text{rk}(\mathcal{E}).\]
The notion of stability then exists for any compact complex manifold:
\vspace{.1in}
\newline
{\em A torsion-free coherent sheaf $\mathcal{E}$ on $M$ is {\em stable} 
if and only if for every coherent subsheaf $\mathcal{S} \subset \mathcal{E}$ 
with $0 < \text{rk}(\mathcal{S}) < \text{rk}(\mathcal{E})$, we have
$\mu(\mathcal{S}) < \mu(\mathcal{E})$.}
\begin{remark} With this definition of stability, many of the 
properties from the K\"{a}hler case hold. For example, 
all line bundles are stable whereas decomposable bundles are always unstable.
In addition, for rank two vector bundles on a surface, 
it is sufficient to verify stability with respect to line bundles;
in particular, if such a bundle is non-filtrable, then it is 
automatically stable. 
Finally, if a vector bundle $E$ is stable, then it is simple, that is,
$h^0(M,End(E)) = 1$. 
\end{remark}

\begin{example}\label{degree of line bundles}
Let $X$ be the classical Hopf manifold corresponding the diagonal 
$(\mu, \dots, \mu)$.
In this case, the degree of line bundles can be 
computed explicitly (for details in the case of surfaces, see \cite{LT,T});
it is determined by a map from ${\rm Pic}(X)$ to 
the reals, denoted $\deg : {\rm Pic}(X) \rightarrow \mathbb{R}$, 
of the form $z \mapsto C\ln|z|$, where $C$ is a real constant.
We define the degree of the line bundle $L_a$, for $a \in \mathbb{C}^\ast$, as 
\[ \deg{L_a} =  \ln|a|/\ln|\mu|, \]
with the normalisation chosen so that 
$\deg{\pi^\ast(\mathcal{O}_{\mathbb{P}^n}(m))} = \deg{L_{\mu^n}} = m$. 
\end{example}

\begin{example}
For a generic Hopf manifold $X$ given by the diagonal $(\mu_1, \dots, \mu_n)$,
$0 < |\mu_1| \leq \dots \leq |\mu_n| < 1$,
we define the degree of the line bundle $L_a$, for $a \in \mathbb{C}^\ast$, as
\[\deg{L_a} = - \ln|a|\]
so that, given any positive integer $m$, we have $\deg{L_{\mu_i^m}} \geq 0$ 
for all $i = 1, \dots, n$. 
In fact, we shall see in section \ref{moduli spaces}
that if we were to define the degree using instead a positive constant $C$, 
then every filtrable rank-2 vector bundle on a generic Hopf surface 
would be unstable.
\end{example}

\section{Holomorphic vector bundles on Hopf surfaces}

\subsection{Topologically trivial holomorphic vector bundles}

We begin by considering topologically trivial holomorphic
vector bundles of arbitrary rank on diagonal Hopf surfaces,  
proving that they possess filtrations by vector bundles.
 
\begin{proposition}\label{non-simple}
On a Hopf surface $X$,
topologically trivial holomorphic vector bundles of rank greater than 1
are not simple. 
\end{proposition}
\begin{proof}  
We prove the proposition by contradiction. 
Consider a holomorphic vector bundle $E$ on $X$ with $c_1(E) = c_2(E) = 0$
and assume that it is simple, implying that $h^0(X;\ad(E)) = 0$.
Recall that the canonical bundle of $X$ is given by $\mathcal{O}_X(-D)$,
where $D$ is the effective divisor $T_1 + T_2$. The inclusion 
$K_X = \mathcal{O}_X(-D) \subset \mathcal{O}_X$ then gives
\[ h^2(X;\ad(E)) = h^0(X;\ad(E) \otimes K_X) = h^0(X;\ad(E)) = 0. \]
Consequently, since $\chi(E) = 0$, we have
$h^1(X;\ad(E) \otimes K_X) = h^1(X;\ad(E)) = 0$.
Inserting this into the long exact sequence on cohomology associated to the exact 
sequence 
\[0 \rightarrow \ad(E) \otimes K_X \rightarrow \ad(E) \rightarrow \ad(E)|_D \rightarrow 0,\]
we obtain $h^0(D;\ad(E)|_D) = 0$,
which contradicts the fact that
$h^0(D;\ad(E)|_D) \geq h^0(T_i;\ad(E)|_{T_i}) \geq 1$.
\end{proof}

Given that we are considering bundles on a surface, 
we then obtain the following.

\begin{corollary}\label{topologically trivial bundles}
Any topologically trivial holomorphic vector bundle on a Hopf surface $X$
possesses a filtration by vector bundles.
\hfill \qedsymbol
\end{corollary}

Holomorphic vector bundles on Hopf manifolds 
that admit filtrations by vector bundles have been studied by
Mall \cite{Mall-filtrable,Mall-rank-2}.
In particular, he shows that they are the only bundles that can
be constructed using factors of automorphy.
Rank-2 vector bundles on a generic Hopf surface $X$ can be classified as follows.

\begin{proposition}[Mall]\label{classification - trivial Chern class}
Let $E$ be an extension of line bundles on $X$.
Then, there exists an exact sequence
\[ 0 \longrightarrow L_a \longrightarrow E \longrightarrow L_b \longrightarrow 0,\]
with $a,b \in \mathbb{C}^\ast$. 
We have the following possibilities:
 
(i) If $a = b \mu_1^{m_1} \mu_2^{m_2}$, for non-negative integers $m_1$ and $m_2$, 
then $E = L_a \oplus L_b$ or $E$ is the {\em unique} non-trivial extension
$0 \longrightarrow L_a \longrightarrow E \longrightarrow L_b \longrightarrow 0$.

(ii) If $a b^{-1} \neq \mu_1^{m_1} \mu_2^{m_2}$ for all integers 
$m_1,m_2 \geq 0$, then $E = L_a \oplus L_b$.
\end{proposition}

\begin{remark}\label{factors of automorphy}
The rank-2 vector bundles described in Proposition 
\ref{classification - trivial Chern class} 
are given by the factors of automorphy:
\begin{equation}\label{factor} 
 \begin{pmatrix}
        a & \epsilon z_1^{m_1}z_2^{m_2}\\
        0 & b
       \end{pmatrix}, 
\end{equation}
where $\epsilon = 0$ if the bundle is decomposable and $\epsilon = 1$ otherwise 
\cite{Mall-filtrable}.
\end{remark}

\subsection{Constructing rank-2 vector bundles}
\label{constructing bundles}

There are three standard ways for constructing rank-2 vector bundles on 
a surface $X$.

{\em (i) Double covers.}
One method for constructing rank-$2$ vector bundles on a surface $X$ is the following.
Find a smooth double cover $\varphi: Y \rightarrow X$ of $X$. 
Then, for any line bundle $L$ on $Y$, the direct image $\varphi_\ast L$ is a rank-$2$
vector bundle.

{\em (ii) Serre construction.}
This method consists in finding locally free extensions of the form
\[ 0 \longrightarrow L \longrightarrow E \longrightarrow 
L' \otimes I_Z \longrightarrow 0, \]
where $L$ and $L'$ are line bundles on $X$ and $I_Z$ is the ideal sheaf $Z$ 
of a finite set of points (counting multiplicity) on $X$ that may be empty.
Bundles thus obtained are clearly all {\em filtrable}.

Note that if $X$ is a Hopf surface, then $c_2(E) = l(Z)$, that is, 
it is equal to the number of points in $Z$ (counting multiplicity).
In addition, extensions of this type exist for any choice of line bundles
$L$ and $L'$ and points on $X$.

\begin{remark*}
Every rank-2 vector bundle on an algebraic surface is filtrable and
can be realised via the Serre construction. On non-algebraic
surfaces, there exist, however, {\em non-filtrable} bundles. In fact, we shall see that
rank-2 vector bundles on Hopf surfaces are generically non-filtrable.
\end{remark*}

{\em (iii) Elementary modifications.}
Start with a rank-2 vector bundle $E$, an effective divisor $D$ on $X$,
and a line bundle $\lambda$ on $D$ such that the restriction of $E$ to $D$ admits
a projection $p: E|_D \rightarrow \lambda$. Let $i:D \rightarrow X$ denote the
natural inclusion. Then, if we also denote by $p$ the induced projection
$E \rightarrow i_\ast \lambda$ on $X$, where $i_\ast \lambda$ is now a torsion 
sheaf supported on $D$, we have the following exact sequence on $X$:
\[ 0 \longrightarrow \bar{E} \longrightarrow E \stackrel{P}{\longrightarrow}
i_\ast \lambda \longrightarrow 0, \]
where $\bar{E} := \ker{p}$ is a rank-2 locally free sheaf; it is called the 
{\em elementary modification of $E$ induced by $p$}.   
Note that $\bar{E}$ is isomorphic to $E$ away from $D$. 
In addition, $\det(\bar{E}) = \mathcal{O}_X(-D) \otimes \det{E}$ and, on a 
Hopf surface, we have $c_2(\bar{E}) = c_1(L) + c_2(E)$.
Other properties of elementary
modifications can be found, for example, in \cite{F}

\begin{remark*}
We shall see that on elliptically fibred Hopf surfaces, all rank-2 vector 
bundles can be obtained this way, but that these methods only produce filtrable
vector bundles on generic Hopf surfaces. In the latter case, the only
method we know so far for constructing non-filtrable vector bundles
is to choose vector bundles on the open cover $X_1,X_2$ of $X$ that
are isomorphic on the overlap and gluing them to obtain a vector
bundle on $X$. Unfortunately, this method makes their classification very difficult.
\end{remark*}

\subsection{Notation and terminology}
To study bundles on a Hopf surface $X$, one of our main tools is 
restriction to one of its elliptic curves $T = \mathbb{C}^\ast / \mu$. 
It is important to point out that the restriction of any vector bundle $E$ on $X$ 
to $T$ is generically semistable, given by an extension of line bundles of degree zero;
if these line bundles correspond to the factors of automorphy $a$ and $b$ in 
$\mathbb{C}^\ast / \mu$, we say that $E$ has {\em splitting type} $(a,b)$  on $T$.
In fact, the restriction of a vector bundle
is unstable on at most an isolated set of curves on $X$. 
If the restriction of $E$ to $T$ is unstable, we say that the vector bundle has
a {\em jump} over $T$.

Consider a rank-2 vector bundle $E$ on $X$ 
with determinant $\delta$ that has a jump of multiplicity $m$ over the curve $T$.
The restriction of $E$ to $T$ is then of the 
form $\lambda \oplus (\lambda^\ast \otimes \delta_{x_0})$, for some 
$\lambda \in \text{Pic}^{-h}(T)$, $h>0$;
the integer $h$ is called the {\it height} of the jump at $T$.
Moreover, up to a multiple of the identity, 
there is a {\em unique} surjection $E|_T \rightarrow \lambda$,
which defines a canonical elementary modification of $E$ that we denote $\bar{E}$; 
this elementary modification is called {\em allowable} \cite{F}.
Therefore, we can associate to $E$ a finite sequence
$\{ \bar{E}_1, \bar{E}_2, \dots , \bar{E}_l \}$
of allowable elementary modifications such that $\bar{E}_l$ is
the only element of the sequence that does not have a jump at $T$. 
The integer $l$ is called
the {\it length} of the jump at $T$. 

\begin{note}
Note that if a vector bundle $E$ jumps over the curve $T$ with mulplicity $m$,
then $m = \sum_{i=1}^{l-1} h_i$, where $h_0 = h$
is the height of $E$ and $h_i$ is the height of
$\bar{E}_i$, $i = 0, \dots, l-1$.
Moreover, if $E$ has $k$ jumps of multiplicity $m_1, \dots, m_k$, respectively,
then $\sum_{i=1}^{k} m_i = c_2(E)$.
For a detailed description of jumps, we refer the reader to 
\cite{Moraru,Brinzanescu-Moraru2}.
\end{note}

\subsection{Moduli spaces}
\label{moduli spaces}
For a fixed line bundle $L_\delta$ on $X$, let
$\mathcal{M}_{\delta,c_2}$ be the moduli space of stable holomorphic rank-2 vector bundles 
with determinant $L_\delta$ and second Chern class $c_2$.
Referring to Proposition \ref{non-simple}, 
the moduli space $\mathcal{M}_{\delta,c_2}$ is non-empty if only if $c_2 > 0$, 
in which case, a similar computation shows that it is a complex manifold of dimension $4c_2$.
 
Stable bundles on elliptic Hopf surfaces have been described in detail 
in \cite{Moraru,Brinzanescu-Moraru2}. Consequently, we focus our presentation
on the generic case, 
briefly stating the results for classical Hopf surfaces that will be 
needed in section \ref{Poisson structures}.
Recall that non-filtrable bundles are automatically stable. 
We therefore only determine stability conditions for the filtrable ones.

\subsubsection{Filtrable rank-2 vector bundles}
\label{filtrable bundles}

Consider a filtrable rank-2 vector bundle $E$ on $X$ with $c_2(E) = c_2 > 0$.
It can therefore be expressed as an extension of the form:
\[ 0 \longrightarrow L_a \longrightarrow E \longrightarrow L_b \otimes I_Z 
\longrightarrow 0,\]
where $Z$ is a set of $c_2$ points (counting multiplicity). 
On a classical Hopf surface, every filtrable bundle can also be
constructed by starting with a bundle with trivial Chern class and adding jumps to it 
to obtain a bundle with the desired second Chern class and ideal sheaf $I_Z$.
This is done by performing elementary modifications, using the fibres of the elliptic
fibration $\pi$ that contain the points of $Z$. In particular, this is always possible
because every point on the surface lies on an elliptic curve.
The advantage of thinking of filtrable bundles this way is that it 
enables us to completely classify them.
Unfortunatly, since generic Hopf surfaces possess only  two curves,
many filtrable vector bundles cannot be constructed this way. 
We nevertheless have a good description of stable filtrable rank-2 vector 
bundles; necessary and sufficient conditions for stability can be
stated as follows. 

\begin{theorem}\label{stable bundles - generic surfaces}
Consider a filtrable rank-2 vector bundle $E$ on a generic Hopf surface 
that has determinant $\delta$ and jumps on $T_1$ and $T_2$ of lengths $l_1$ and $l_2$,
respectively.
Suppose that $L_a$ is one of the line bundles of maximal degree mapping
into $E$.
Then, $E$ is stable if and only if 
\begin{equation}\label{condition a} 
a^2 = \delta \mu_1^{-l_1 -k_1} \mu_2^{-l_2 - k_2} 
\end{equation}

\noindent
for two non-negative integers $k_1$ and $k_2$ that are not both zero, or 
\[ a \in D_{l_1,l_2} := \left\{ \alpha \in \mathbb{C}^\ast \ : \ 
|\delta|^{1/2} < |\alpha| < |\delta \mu_1^{-2l_1}\mu_2^{-2l_2}|^{1/2} \right\}. \]
In particular, this implies that $l_1 + l_2 > 0$ unless $a$ satisfies 
equation \eqref{condition a}.
\end{theorem}
\begin{proof}
Recall that $E$ is stable if and only if each of its destabilising bundles 
has degree less than $\deg(L_\delta)/2$.
Consider the rank-2 vector bundle $\bar{E}$ obtained by performing
$l$ elementary modifications to remove the jumps of $E$.
Then, $\det \bar{E} = L_{\delta \mu_1^{-l_1}\mu_2^{-l_2}}$.
Note that $E$ and $\bar{E}$ have the same destabilising bundles.
Indeed, $\bar{E}$ is obtained by taking consecutive elementary modifications
determined by exact sequences of the form:
\[ 0 \rightarrow \bar{E}_{i+1} \rightarrow \bar{E}_i \rightarrow i_\ast \lambda \rightarrow 0,\]
where $\lambda$ is a line bundle of negative degree on $T_1$ or $T_2$, 
for $0 \leq i \leq l_s - 1$, $s = 1,2$.
Thus, $h^0(X, L_{c^{-1}} \otimes i_\ast\lambda) = 0$ for all line bundles $L_c$,
so that 
\[ h^0(X,L_{c^{-1}} \otimes \bar{E}_{i+1}) = h^0(X,L_{c^{-1}} \otimes \bar{E}_{i}).\] 

\noindent
This implies, in particular, that $\bar{E}$ is an extension
of $ L_{a^{-1}\delta \mu_1^{-l_1}\mu_2^{-l_2}} \otimes I_{Z'}$ by $L_a$,
where $Z'$ is a (possibly empty) set of points that do not lie on $T_1$ or $T_2$. 

We therefore have to determine the line bundles that map non-trivially into $\bar{E}$.
Let $L_c$ be such a line bundle.
We first assume that $Z'$ is empty.
Suppose that $\bar{E}$ decomposes as $L_a \oplus L_{a^{-1}\delta \mu_1^{-l_1}\mu_2^{-l_2}}$; 
then $L_a$ and $L_{a^{-1}\delta \mu_1^{-l_1}\mu_2^{-l_2}}$ are the line bundles
of maximal degree mapping into $\bar{E}$.
To ensure the stability of $E$, both must have degree strictly smaller than $\deg (L_\delta) /2$.
Referring to the definition of degree given in section \ref{degree},
this is equivalent to $a$ being an element of $D_{l_1,l_2}$.

If $\bar{E}$ is instead indecomposable,
then there exist $a \in \mathbb{C}^\ast$ and non-negative integers $m_1,m_2$ 
such that $\bar{E}$ is given by the non-trivial extension of 
$L_{a \mu_1^{-m_1} \mu_2^{-m_2}}$ by $L_a$ 
(see Proposition \ref{classification - trivial Chern class} (i)).
Note that any line bundle mapping into $\bar{E}$ must also map to $L_a$.
Indeed, if $h^0(X,L_{c^{-1}}\bar{E}) \neq 0$, then $h^0(X,L_{c^{-1}a})$ and
$h^0(X,L_{c^{-1}a\mu_1^{-m_1}\mu_2^{-m_2}})$ cannot both zero;
hence, since $h^0(X,L_{c^{-1}a\mu_1^{-m_1}\mu_2^{-m_2}}) \leq h^0(X,L_{c^{-1}a})$, 
we have $h^0(X,L_{c^{-1}a}) \neq 0$
and $c^{-1}a = \mu_1^{k_1}\mu_2^{k_2}$ for some integers $k_1,k_2 \geq 0$.
Consequently, $L_a$ is the (unique) line bundle of maximal degree mapping into $\bar{E}$.
Suppose that it has degree strictly smaller than $\deg (L_\delta) /2$.
Then, since 
\[ a^2 \mu_1^{-m_1}\mu_2^{-m_2} = \delta \mu_1^{-l_1}\mu_2^{-l_2}\]
and $|\mu_1^{-m_1}\mu_2^{-m_2}| > 1$, it follows that $a$ is an element of
$D_{l_1,l_2}$. 

Let us now assume that $Z'$ is not empty. Note that $\bar{E}$ is semistable
on both $T_1$ and $T_2$. To simplify the notation, let us set 
$b = a^{-1}\delta \mu_1^{-l_1}\mu_2^{-l_2}$.
The set of line bundles mapping into $\bar{E}$ is therefore contained in 
\[ \left\{ L_c \  :  \ \mbox{$c = a\mu_1^{-k_1}\mu_2^{-k_2}$ \ \   or \ \   
$c = b \mu_1^{-k_1}\mu_2^{-k_2}$ for integers $k_1,k_2 \geq 0$} \right\}.\]
Let $L_c$ be another destabilising line bundle of $\bar{E}$
so that $E / L_c$ is torsion free.
Then $c = b \mu_1^{-k_1}\mu_2^{-k_2}$ for integers $k_1,k_2 \geq 0$.
If at least one of the integers is non-zero, then we must have $c = a$,
otherwise the quotient $E / L_c$ is not torsion-free 
(recall that on a generic Hopf surface, there are no non-trivial relations
between $\mu_1$ and $\mu_2$).
Consequently, $L_a$ is the (only) line bundle of maximal degree mapping into
$E$; since $a = a^{-1}\delta \mu_1^{-l_1-k_1}\mu_2^{-l_2-k_2}$
with $k_1$ or $k_2$ non-zero, we have $|a| > |\delta|^{1/2}$, implying
that $E$ is stable. Finally, if $k_1 = k_2 = 0$, then the second destabilising 
line bundle of $E$ is $L_{a^{-1}\delta \mu_1^{-l_1}\mu_2^{-l_2}}$
so that $E$ is stable if and only if $a \in D_{l_1,l_2}$.
\end{proof}

\begin{remark}\label{stable bundles - classical}
On a classical Hopf surface $X$, determined by a
diagonal $(\mu,\mu)$, vector bundles can have jumps on elliptic curves other 
than $T_1$ and $T_2$. 
Let $E$ be a filtrable rank-2 vector bundle on $X$
that has determinant $L_\delta$ and $k$ jumps of 
lengths $l_1, \dots, l_k$, respectively. Set $l = l_1 + \dots + l_k$. 
Using the notation of Theorem \ref{stable bundles - generic surfaces},
an extension of $L_{a^{-1} \delta} \otimes I_Z$ by $L_a$
is stable if and only if $a \in D_l$, where
\[ D_l := \left\{ \alpha \in \mathbb{C}^\ast \ : \ 
|\delta|^{1/2} < |\alpha| < |\delta \mu^{-2l}|^{1/2} \right\}. \]
This follows from Lemma 4.5 and Corollary 4.6 of \cite{Moraru}
(with $\mu = \lambda^{-1}$ because classical Hopf surfaces were defined
by a complex number $\lambda$ with $|\lambda| > 1$ in \cite{Moraru}).
\end{remark}

\begin{remark}
The domains $D_{l_1,l_2}$ and $D_l$ defined in Theorem 
\ref{stable bundles - generic surfaces} and Remark \ref{stable bundles - classical}
are independent of the definition of degree, up to multiplication by a positive
constant. Otherwise, one readily verifies that these domains would in fact be
empty. 
\end{remark}

We now describe stable 
filtrable bundles with $c_2 = 1$ and fixed determinant $\delta$
that have a jump on $T_1$ or $T_2$.
Without loss of generality, we assume it to be $T_1$.
Note that a similar analysis can be carried out for bundles with $c_2 > 1$.

\begin{proposition}\label{stable filtrable c_2 = 1}
Let $E$ be a stable filtrable rank-2 vector bundle on $X$
with determinant $\delta$ and a jump of multiplicity $1$ on $T_1$. 
Then, $E$ is uniquely determined by a triple $(a,\lambda,p)$ such that
\[ (a,\lambda) \in D_{1,0} \times {\rm Pic}^1(T_1)\]
and $p$ is a projection from $L_{a^{-1}\delta} \oplus L_{a\mu_1}$
to $\lambda$ on $T_1$ that is unique up to isomorphism, unless 
$a^2 = \delta \mu_1^{m_1-1}\mu_2^{m_2}$ with 
$m_1 \geq 1$ and $m_2 > 0$,
in which case it is an element of the projective line
$\mathbb{P}^1(H^0(T_1, Hom(L_{a^{-1}\delta} \oplus L_{a\mu_1}, \lambda)))$.

Note that $L_a$ corresponds to one of the destabilising line bundles 
of $E$ and that $\lambda$ is such that the restriction of $E$ to $T_1$ splits as 
$\lambda \oplus (\lambda^{-1} \otimes \delta)$.
\end{proposition}
\begin{proof}
Let $E$ be a stable filtrable rank-2 vector bundle on $X$
with a jump of multiplicity $1$ on $T_1$. 
Such a bundle $E$ is then given by an extension of the
form 
\[ 0 \rightarrow L_a \rightarrow E \rightarrow L_{a^{-1}\delta} \otimes I_p, \]
where $a \in D_{1,0}$ and $p$ is a point on $T_1$.
Moreover, the allowable elementary modification $\bar{E}$ of $E$ 
is an extension of $L_{a^{-1}\delta \mu_1^{-1}}$ by $L_a$
that splits unless $a^2 = \delta \mu_1^{m_1-1}\mu_2^{m_2}$ with 
$m_1 \geq 1$ and $m_2 > 0$ (note that if $m_2 = 0$, then $a \notin D_{1,0}$).
In the latter case, $\bar{E}$ is given by a factor of automorphy of
the form \eqref{factor} with $\epsilon = 0$ or $1$.

Suppose that the splitting type of $E$ on $T_1$ is 
$\lambda \oplus (\lambda^{-1} \otimes \delta)$ 
for some $\lambda \in {\rm Pic}^1(T_1)$.
Then, $E$ can be recovered from $\bar{E}$ by using the line bundle $\lambda$ 
to introduce a jump to $\bar{E} \otimes L_{\mu_1}$ on $T_1$. 
Given a fixed 
choice of line bundle $\lambda$ in ${\rm Pic}^1(T_1)$, we therefore have to determine 
which projections $p$ from 
$\bar{E} \otimes L_{\mu_1}|_{T_1}$ to $\lambda$
induce isomorphic elementary modifications.
Note that since no element $a$ in $D_{1,0}$ is such that 
$a^2 \equiv \delta \mod \mu_1$, such projections exist for all extension
of $L_{a^{-1}\delta}$ by $L_{a\mu_1}$.
In addition, any two such projections differ an
element of $Aut(\bar{E} \otimes L_{\mu_1}|_{T_1}$)
because $\deg \lambda = 1$; 
it is therefore sufficient to
find out which automorphisms of $\bar{E} \otimes L_{\mu_1}|_{T_1}$ 
extend to automorphisms of $\bar{E} \otimes L_{\mu_1}$ on $X$.

If $\bar{E} \otimes L_{\mu_1}$ is decomposable, then 
$Aut(\bar{E} \otimes L_{\mu_1}|_{T_1}) = Aut(\bar{E} \otimes L_{\mu_1})$;
in this case, there is a unique way, up to isomorphism, of introducing the jump
so that $E$ is uniquely determined by the pair $(a,\lambda)$.
Let us now assume that $\bar{E} \otimes L_{\mu_1}$ is a non-trivial
extension of $L_{a^{-1}\delta}$ by $L_{a\mu_1}$, where $a^2 = \delta \mu_1^{m_1-1}\mu_2^{m_2}$ with 
$m_1 \geq 1$ and $m_2 > 0$. 
In this case, the restriction of $\bar{E} \otimes L_{\mu_1}$ to $T_1$ is decomposable
and the only elements of $Aut(\bar{E} \otimes L_{\mu_1}|_{T_1})$
that extend to $Aut(\bar{E} \otimes L_{\mu_1})$ are multiples of the identity.  
Consequently, the vector bundle $E$ is determined by a triple $(a,\lambda,p)$
such that $a \in D_{1,0}$, $\lambda \in {\rm Pic}^1(T_1)$, and 
$p$ is a projection in 
$\mathbb{P}^1(H^0(T_1, Hom(L_{a^{-1}\delta} \oplus L_{a\mu_1}, \lambda)))$.
\end{proof}

\begin{remark}\label{stable filtrable c_2 = 1 - classical}
Note that if $|\mu_1| = |\mu_2|$, then the description simplifies.
In this case, the allowable elementary modification 
$\bar{E}$ always decomposes as
$L_a \oplus L_{a^{-1}\delta\mu_1^{-1}}$, otherwise $a \notin D_{1,0}$; 
stable filtrable bundles with $c_2 = 1$, determinant $\delta$,
and a jump on $T_1$ are therefore in one-to-one correspondence with the pairs 
in $D_{1,0} \times {\rm Pic}^1(T_1)$.
In addition, if we consider such bundles
on a classical Hopf surface $X = \mathbb{C}^2_\ast/(\mu,\mu)$,
then we obtain the same description, except that
the bundles can now have a jump over any fibre of the elliptic fibration. 
For a fixed fibre $T$, these bundles are parametrised by
$D_1 \times {\rm Pic}^1(T)$.
\end{remark}

As direct consequence of the above discussion, we have:

\begin{corollary}
Consider the moduli space $\mathcal{M}_{\delta,1}$ of stable rank-2 vector bundles 
on a Hopf surface $X$ with determinant $\delta$ and second Chern class $1$.
Let $F\mathcal{M}_{\delta,1}$ be the subset of $\mathcal{M}_{\delta,1}$
consisting of filtrable bundles.
Then every component of $F\mathcal{M}_{\delta,1}$ has codimension at least 
one in $\mathcal{M}_{\delta,1}$.
\hfill \qedsymbol
\end{corollary}

\begin{example}\label{magnetic monopoles}
We end this section with an application of the above analysis 
to magnetic monopoles on solid tori. These can be seen to correspond
to $S^1$-invariant instantons on Hopf surfaces 
of the form
\[ \mathbb{C}^2_\ast / (\mu,|\mu|), \]
where $\mu$ is a complex number with $|\mu| < 1$ \cite{Braam1,Braam2}.
Monopoles on solid tori can therefore also be identified, via 
the Hitchin-Kobayashi correspondence \cite{LT}, 
with $\mathbb{C}^\ast$-equivariant stable holomorphic bundles on Hopf surfaces.

Moduli spaces of $\mathbb{C}^\ast$-equivariant stable holomorphic bundles on 
Hopf surfaces were first studied by Braam and Hurtubise \cite{B-H}
in the classical case.
They showed that the moduli spaces $\mathcal{M}(m,k)$
of monopoles of mass $m$ and charge $k$ 
are complex manifolds of dimension $2k$ consisting of certain
stable filtrable rank-2 vector bundles with $c_2 = mk$ and a jump of height $k$ and
length $m$ over $T_1$. 
In particular, they showed that $\mathcal{M}(m,1)$
is isomorphic $D_m \times {\rm Pic}^1(T_1)$, by classifying $\mathbb{C}^\ast$-equivariant
bundles on $X_1$ and $X_2$ and determining how one can glue them on the overlap
to obtain distinct monopoles.
The method presented in Proposition \ref{stable filtrable c_2 = 1} offers, however, a more
invariant way of approaching the problem. In fact, a similar analysis 
shows that for $k>1$, the moduli space $\mathcal{M}(m,k)$
consists of triples $(a,\lambda,p)$, where 
$a \in D_m$, $\lambda \in {\rm Pic}^k(T_1)$, 
and $p$ is a projection on $T_1$ from $L_{a\mu^m} \oplus L_{a^{-1}\delta}$ to $\lambda$
or possibly from the non-trivial extension of 
$L_{a^{-1}\delta}$ by $L_{a\mu^m}$ to $\lambda$ when $a^2 \equiv \delta^{-1}\mod \mu$.
The space of such projections is in this case $(2k-2)$-dimensional.
\end{example}

\subsubsection{Non-Filtrable bundles}
\label{non-filtrable bundles}

In this section, we turn to the problem of constructing non-filtrable 
holomorphic vector bundles.
We begin by noting that there exist many non-filtrable bundles on Hopf surfaces. 
For example, the generic elements of $\mathcal{M}_{\delta,1}$ are non-filtrable
because the set of filtrable bundles in $\mathcal{M}_{\delta,1}$ 
has codimension at least one 
(see sections \ref{constructing bundles} and \ref{filtrable bundles}).  
Such vector bundles can then be used to construct non-filtrable
vector bundles of arbitrary second Chern class by performing elementary modifications
to add (or increase) jumps. 

On elliptically fibred Hopf surfaces, non-filtrable holomorphic 
vector bundles are by now well-understood.
In fact, one can show that they can all be constructed by using double
covers and elementary modifications
(for details, see 
\cite{Moraru,Brinzanescu-Moraru,Brinzanescu-Moraru1}).
For the convenience of the reader, we briefly recall how this is done in the case 
of a classical Hopf surface.

Let $X$ be a classical Hopf surface.
Referring to section \ref{preliminaries}, it admits an elliptic fibration
$\pi: X \rightarrow \mathbb{P}^1$, with fibre an elliptic curve $T$
and relative Jacobian $J(X)=\mathbb{P}^1\times T^*$.  
Consider a vector bundle $E$ on $X$. One of the main tools for studying this bundle
is restriction to the fibres of the fibration $\pi$.
In particular, there exists a divisor $S_E$ in the relative Jacobian of $X$, 
called the {\em spectral curve} or {\em cover} of the bundle, 
that encodes the isomorphism
class of the bundle $E$ over each fibre of $\pi$.
Note that the self-intersection of this divisor is equal to a multiple of the second
Chern class of the bundle.

\begin{example*}
Let $E$ be a rank-2 vector bundle on $X$ with determinant $L_\delta$ and
second Chern class $c_2$. Then it has a spectral curve of the form
\[ S_E := \left( \sum_{i=1}^k \{ x_i \} \times T^\ast \right) + C_E, \]
where $C_E$ is a bisection of $J(X)$ and $x_1, \cdots, x_k$ are points in 
$\mathbb{P}^1$ corresponding to the jumps of $E$. 
Away from the jumps, the pair of points 
$(\lambda_0,\lambda_0^{-1} \otimes L_\delta)$
on $C_E$ above $x_0 \in \mathbb{P}^1$ gives the splitting
type of $E$ on the fibre $T_{x_0} = \pi^{-1}(x_0)$.
In this case, 
\[ S_E \cdot S_E = 4c_2. \]
Moreover, if the bisection $C_E$ is smooth, then it is a double
cover of $\mathbb{P}^1$ of genus $(2c_2 - 2k - 1)$.
\end{example*}

Let us now describe non-filtrable vector bundles on $X$. 
We have seen that we can associate to any rank-2 vector bundle $E$ on $X$ a bisection 
$C_E \subset J(X)$. 
The vector bundle $E$ is then filtrable if and only if its 
bisection $C_E$ is reducible.
Conversely, given any bisection $C$ of $J(X)$, one can associate to it at least one 
rank-2 vector bundle on $X$. 
This implies, in particular, that 
non-filtrable bundles exist on classical Hopf surfaces because 
irreducible bisections exist in $J(X)$. 
The bundles determined by $C$ are constructed as follows. Consider the double cover 
\[ \varphi: Y := X \times_{\mathbb{P}^1} C \rightarrow X; \] 
then, for any line bundle $L$ on $Y$, 
the pushdown $\varphi_\ast L$ is a rank-2 vector bundle on $X$. 
In fact, for a certain class of line bundles on $Y$, the resulting
rank-2 vector bundles will have spectral cover $C$.
For example, if the bisection is smooth, then one can show that
the bundles that correspond to it are parametrised by
the abelian variety $Jac(C)$ (see \cite{Moraru} for precise statements).

\begin{note}
The spectral construction applies, in fact, to any elliptic 
fibration; it has been used by many authors to study bundles on elliptic fibrations
(see for example \cite{F,FM,FMW,D,T}).
\end{note}

Any holomorphic rank-2 vector bundle on a classical Hopf surface
can therefore be constructed by using a double cover and 
elementary modifications (to add jumps).
This is, however, certainly not the case for generic Hopf surfaces.

\begin{proposition}
On a generic Hopf surface $X = \mathbb{C}^2_\ast/(\mu_1,\mu_2)$, only filtrable
vector bundles can be constructed by using double covers.
\end{proposition}
\begin{proof}
Consider a double cover $\varphi:Y \rightarrow X$ of $X$
determined by the line bundle $\mathcal{L}$ on $X$.
Then, for any line bundle $\mathcal{M}$ on $X$, we have
\[ c_2(\varphi_\ast(\mathcal{M}) = \frac{1}{2}(c_1^2({\rm Nm}\mathcal{M})) - 
\varphi_\ast(c_1^2(\mathcal{M}))-\varphi_\ast(c_1(\mathcal{M})).c_1({\mathcal L})) \]
(see \cite{Brinzanescu}). Therefore, since $h^2(X,\mathbb{Z}) = 0$, this reduces to
\begin{equation}\label{c2}
c_2(\varphi_\ast(\mathcal{M}) = -\varphi_\ast(c_1^2(\mathcal{M})). 
\end{equation}
Without loss of generality, we can assume that $Y$ is a smooth surface
(otherwise, take its normalisation). Referring to Lemma \ref{double covers},
$Y$ is then either a generic Hopf surface or one of the non-primary Hopf surfaces
$\Theta^\ast_{-2}/(\beta,\mu_1)$, with $\beta^2 = \mu_1$, 
described in section \ref{non-primary-Hopf}.
This means, in particular, that the first Chern class of any line bundle $L$ on
$X$ is torsion. Consequently, by \eqref{c2}, we have $c_2(\varphi_\ast(L)) = 0$.
The rank-2 vector bundle $\varphi_\ast L$ is thus filtrable by 
Corollary \ref{topologically trivial bundles}.
\end{proof}

\subsection{Poisson structures}
\label{Poisson structures}

A (holomorphic) Poisson structure on a complex surface is given by a global 
section of its anticanonical bundle \cite{Bottacin-Poisson}.
Consequently, any Hopf surface $X$ admits a Poisson structure because its
anticanonical bundle $K_X^{-1}$ is given by the effective divisor $D := T_1 + T_2$. 
Fix a Poisson structure $s \in H^0(X,K_X^{-1})$ on $X$. 
A Poisson structure $\theta = \theta_s \in 
H^0(\mathcal{M},\otimes^2T\mathcal{M})$ on the moduli space 
$\mathcal{M} := \mathcal{M}_{\delta,c_2}$ is then defined as follows: for any bundle
$E \in \mathcal{M}$, 
$\theta(E) : T^\ast_E\mathcal{M} \times T^\ast_E\mathcal{M}
\longrightarrow \mathbb{C}$ 
is the composition
\[ \begin{array}{lr}
\theta(E): H^1(X,\ad(E) \otimes K_{X}) \times
H^1(X,\ad(E) \otimes K_{X}) 
\stackrel{\circ}{\longrightarrow} \\
\hspace{0.9in} H^2(X,\End(E) \otimes K^2_{X}) 
\stackrel{s}{\longrightarrow}
H^2(X,\End(E) \otimes K_{X}) 
\stackrel{\text{Tr}}{\longrightarrow} \mathbb{C},
\end{array} \]
where the first map is the cup-product of two cohomology classes, the second is 
multiplication by $s$, and the third is the trace map. 

The Poisson structure $s$ is degenerate, its divisor being $D = (s)$.  
Moreover, at any point $E \in \mathcal{M}$, 
\[ \rk \theta(E) = 4c_2  - \dim H^0(D,\ad(E|_D)). \]

\noindent
We see that the rank of the Poisson structure is generically $4c_2 - 2$, 
and ``drops'' at the points of $\mathcal{M}$ corresponding to bundles that 
are not regular over the fibres $T_1$ and $T_2$ 
(for details in the elliptic case, see \cite{Moraru}).

On the set of bundles for which the Poisson structure is maximal,
one can define the following maps:
\[ \begin{array}{rcl}
f_i: \mathcal{M}_{\delta,c_2} & \longrightarrow & 
\text{Pic}^0(T_i)/ i_{\delta} = \mathbb{P}^1 \\
E & \mapsto & (L_{a_i}|_{T_i}, L_{a_i^{-1}\delta}|_{T_i}),
\end{array} \]

\noindent
where $L_{a_i}$ and $L_{a_i^{-1}\delta}$ are the destabilising bundles of $E|_{T_i}$,
and $i_{\delta}$ is the involution of $\text{Pic}^0(T_i)$ given by
$\lambda_0 \mapsto \lambda_0^{-1} \otimes L_\delta$.
The functions $f_1$ and $f_2$ are linearly independent Casimirs;
this can be proven as in the elliptic case \cite{Moraru}.
\vspace{.05in}

On a classical Hopf surface $X$, the Casimirs can be described very explicitly
in terms of the spectral data of the bundle.
Let $E$ be rank-2 vector bundle on $X$ with determinant $L_\delta$
and spectral curve $S_E \subset J(X)$ (see section \ref{non-filtrable bundles}).
In this case, one can associated to $E$ an equivalent divisor 
that is constructed as follows.
Consider the quotient of $J(X) = \mathbb{P}^1 \times T^\ast$ by 
the involution $i_\delta : = id \times i_\delta$, where $id$ is the identity
on $\mathbb{P}^1$
and $i_\delta$ is the involution of $T^\ast$ determined by $L_\delta$ (see above).
Let $\eta: J(X) \rightarrow J(X)/i_\delta = \mathbb{P}^1 \times \mathbb{P}^1$ 
be the canonical map.
By construction, the spectral curve of $E$ is invariant with respect to
the involution $i_\delta$. 
It therefore descends to a divisor on 
$\mathbb{P}^1 \times \mathbb{P}^1$ of the form 
\[ \mathcal{G}_E := \left( \sum_{i=1}^k \{ x_i \} \times \mathbb{P}^1 \right) + \text{Gr}(F_E),\]
where $\text{Gr}(F_E)$ is a the graph of a rational map 
$F_E: \mathbb{P}^1 \rightarrow \mathbb{P}^1$ of degree $(c_2(E) - k)$
such that $\eta^\ast \text{Gr}(F_E) = C_E$. This divisor is called the {\em graph}
of $E$; it is an element of the linear system 
$|\mathcal{O}_{\mathbb{P}^1 \times \mathbb{P}^1}(c_2,1)|$.
Note that the bundle $E$ is filtrable if and only if the map $F_E$ is constant.
\vspace{.05in}

Let $x_1$ and $x_2$ be the points in $\mathbb{P}^1$ such that 
$\pi^{-1}(x_i) = T_i$, for $i=1,2$.
The Casimirs $f_i: \mathcal{M}_{\delta,c_2} \rightarrow \mathbb{P}^1$, $i=1,2$, 
are then given by
$E \mapsto F_E(x_i)$, where $F_E$ is the graph of $E$. 
These functions can be used to describe the symplectic leaves of the Poisson structure; 
the symplectic leaves of maximal rank are labelled as follows:
\[ \mathcal{L}_{C_1,C_1} := 
\{ E \in \mathcal{M}_{\delta,c_2} \ | \ \mbox{$\rk\theta(E) = 4n-2$ and 
$f_i(E) = C_i$ for $i = 1,2$}  \}. \]

\noindent
On the open dense subset of non-filtrable bundle that are regular
on every fibre of $\pi$, the elements of the leaf $\mathcal{L}_{C_1,C_1}$ 
can be identified with pairs of the form $(Gr(F), E)$, where $Gr(F)$ is 
the graph of a rational map 
$F: \mathbb{P}^1 \rightarrow \mathbb{P}^1$ of degree $c_2$ passing through the 
points $(x_1,C_1)$ and $(x_2,C_2)$ in $\mathbb{P}^1 \times \mathbb{P}^1$
and $E$ is a rank-2 vector bundle whose graph is given by $Gr(F)$.
Recall that for such a graph,
the set of all bundles corresponding to it
is parametrised by the Jacobian of a curve of genus $2c_2-1$, 
given by $C = \eta^\ast Gr(F)$, confirming that  
the leaf $\mathcal{L}_{C_1,C_1}$ is $(4c_2-2)$-dimensional.
\vspace{.05in}

It is in fact possible to give an explicit description of all the 
symplectic leaves of the Poisson structure. 
We finish by describing these leaves 
for the moduli space $\mathcal{M}_{\delta,1}$ on a classical Hopf surface.

\begin{example}
Let us first consider stable vector bundles $E$ in $\mathcal{M}_{\delta,1}$ 
where the Poisson structure $\theta$ has maximal rank 2.  
On a classical Hopf surface, $f_1(E) = f_2(E)$ if and only if 
the vector bundle $E$ filtrable 
(because if the map $F_E$ has degree one, then it must be injective).
Let us assume that $C_1 = C_2$; then,
\[ \mathcal{L}_{C_1,C_1} = \{ (x_0,\lambda) \ | \ 
\mbox{$x_0 \in \mathbb{P}^1 \backslash \{x_1,x_2\}$ and $\lambda \in {\rm Pic}^1(T_{x_0})$} \}.\]

\noindent
Indeed, stable filtrable bundles on elliptic Hopf surfaces
with $c_2=1$ are completely determined by the choice of a fibre 
$T_{x_0} = \pi^{-1}(x_0)$, over which they have a jump, and
a pair $(a,\lambda)$ in $D_1 \times {\rm Pic}^1(T_{x_0})$ 
(see Remark \ref{stable filtrable c_2 = 1 - classical}). However,
the Casimirs fix $a$ so that $x_0$ and $\lambda$ are the only free parameters left.
If $C_1 \neq C_2$, then the leaf $\mathcal{L}_{C_1,C_1}$ is given as above.

Finally, at the remaining points of the moduli space, 
the Poisson structure has rank zero.
These points correspond to stable filtrable bundles that have a jump 
of multiplicity one
on either $T_1$ or $T_2$. 
Referring to Remark \ref{stable filtrable c_2 = 1 - classical}, 
these bundles are completely determined
by their restrictions to the curves $T_1$ or $T_2$, which fix the 
pairs $(a,\lambda)$ parameterising them.
\end{example}

\begin{remark}
If $X$ is a non-generic Hopf surface, then the graph map 
\[G: \mathcal{M}_{\delta,c_2} \rightarrow 
|\mathcal{O}_{\mathbb{P}^1 \times \mathbb{P}^1}(c_2,1)|
= \mathbb{P}^{2c_2+1},\] 
which associates to each bundle its graph,
admits a structure of algebraically completely integrable Hamiltonian systems
\cite{Moraru, Brinzanescu-Moraru2}. 
\end{remark}

\section{Holomorphic vector bundles on higher dimensional Hopf manifolds}

We have seen that holomorphic vector bundles on Hopf surfaces are generically 
non-filtrable. In contrast, they are always filtrable on higher dimensional 
Hopf manifolds. More precisely, the following is known in this case \cite{V1,V2}.

\begin{theorem}[Verbitsky]\label{sheaves are filtrable}
Let $X$ be a diagonal Hopf manifold of dimension greater than two.
Then all coherent sheaves on $X$ are filtrable.  
\end{theorem}

\begin{theorem}[Verbitsky]\label{Verbitsky - stable bundles}
Let $\pi: X \rightarrow \mathbb{P}^{n-1}$ be a classical Hopf manifold 
of dimension $n \geq 3$. 
Let $E$ be a stable holomorphic vector bundle $E$ on $X$.
Then $E = L \otimes \pi^\ast E'$, where $L$ is a line bundle on $X$ and 
$E'$ is a stable bundle on $\mathbb{P}^{n-1}$. 
\end{theorem}

These results were first proven for positive principal
elliptic fibrations of dimension greater than two \cite{V1}, examples of which are given by 
classical Hopf manifolds. 
Verbitsky shows that on such an elliptic fibration $\pi:X \rightarrow M$, where $M$ is 
a K\"{a}hler manifold dimension at least two, stable vector bundles are equivariant 
with respect to a torus action. He then uses this to show that
stable bundles on $X$ are of form the $L \otimes \pi^\ast E'$;
when the base $M$ is projective, this implies that all holomorphic
vector bundles are filtrable.
His equivariance argument extends, however, to give filtrability
of bundles on all diagonal Hopf manifolds \cite{V2}.

Note that on a classical Hopf manifold $X$, the filtrability of vector bundles
can easily be seen via the spectral construction presented in section 
\ref{non-filtrable bundles} because bundles on such a manifold are topologically trivial.

\begin{proposition}
Let $E$ be a rank-$r$ vector bundle on $X$. 
Then $E$ is topologically trivial.
Moreover, its spectral cover $S_E$ is an effective
divisor in $J(X) = \mathbb{P}^{n-1} \times T^\ast$ of the form
\begin{equation}\label{spectral cover} 
S_E = \sum_{i = 1}^r \mathbb{P}^{n-1} \times \{ \lambda_i \}, 
\end{equation}
for some points $\lambda_1, \dots, \lambda_r$ in $T^\ast$,
implying that it is filtrable.
\end{proposition}
\begin{proof}
Let $L$ be a line bundle on $X$ such that
$h^0(\pi^{-1}(x),L^\ast \otimes E) = 0$ for generic $x \in \mathbb{P}^{n-1}$.
Then, $R^1\pi_\ast(L^\ast \otimes E)$ is a torsion sheaf on 
$\mathbb{P}^{n-1}$ and $R^i\pi_\ast(L^\ast \otimes E) = 0$ for $i = 0$ and $i \geq 2$.
Suppose that the line bundle $L$ corresponds to the section $X \times \{ \lambda \}$
of $J(X)$. Note that 
\[ \text{Pic}(J(X)) = \text{Pic}(\mathbb{P}^{n-1}) \times \text{Pic}(T^\ast).\] 

\noindent
Recall that the spectral cover $S_E$ of $E$ is an effective divisor on $J(X)$ that is 
an $r$-to-one cover of $\mathbb{P}^{n-1}$. Consequently, 
we have
\[ S_E \sim \sum_{i = 1}^r \mathbb{P}^{n-1} \times \{ \lambda_i \} +  m \{ H \} \times T^\ast\]
for some non-negative integer $m$,
where the $\lambda_i$'s are points in $T^\ast$, $H$ is a hyperplane in $\mathbb{P}^{n-1}$.
Then, $S_E \cdot (X \times \{ \lambda \}) = m \{ H \} \times \{ \lambda \}$
and the support of $R^1\pi_\ast(L^\ast \otimes E)$ is a divisor
on $\mathbb{P}^{n-1}$ equivalent to $mH$. Furthermore, if $h$ is the Poincar\'{e} dual
of $H$ in $H^2(\mathbb{P}^{n-1},\mathbb{Z})$, then
$c_1(R^1\pi_\ast(L^\ast \otimes E)) = mh$. 
Note that $h$ is the positive generator of 
$H^2(\mathbb{P}^{n-1},\mathbb{Z})$.
Given that $ch(E) = r + (-1)^{(n-1)}c_n(E)/(n-1)!$ 
and $td(X) = 1$, 
by Grothendieck-Riemann-Roch, it follows that
\[ ch(R^1\pi_\ast(L^\ast \otimes E)) = \pi_\ast \left( ch(E) \cdot td(X) \right)
\cdot td(\mathbb{P}^n)^{-1} = (-1)^n\frac{c_n(E)}{(n-1)!} h^{n-1}.\]
Consequently, since $n > 2$, we see that $m = c_n(E) = 0$ 
and $S_E$ is of desired the form.
\end{proof}

\noindent
A vector bundle $E$ on a classical Hopf manifold $X$
therefore either admits a filtration by vector bundles, 
which decomposes if all the $\lambda_i$'s appearing in its 
spectral cover $S_E$ \eqref{spectral cover} are distinct, 
or it can also be of the form $L \otimes \pi^\ast E'$
when all the $\lambda_i$'s are equal.
Unfortunatly, this analysis cannot be successfully carried out 
for bundles on generic Hopf manifolds as there is no satisfactory analogue of 
the spectral construction in the generic case (because the are so few divisors
on these manifolds). 
\vspace{.1in}

On elliptically fibred Hopf manifolds of dimension greater than two,
the study of stable vector bundles boils down to the difficult problem 
of classifying stable vector bundles on projective spaces of dimension 
greater than one.
On generic Hopf manifolds of dimension greater than two, 
the question is, however, greatly simplified by the
fact that these manifolds possess few subvarieties.
One can therefore obtain a complete classification of vector bundles 
by studying extensions of sheaves.
For brevity, we restrict our presentation to the case of holomorphic rank-2 vector bundles; 
similar results however hold for bundles of arbitrary rank. 
In particular, we show that there exist stable rank-2 vector bundles on these manifolds.
 
\begin{proposition}\label{classification bundles generic manifold}
Let $X$ be a generic Hopf manifold of dimension $n \geq 3$ given by the
quotient $\mathbb{C}^n_\ast / (\mu_1,\dots,\mu_n)$.
If $E$ be a rank-2 vector bundle on $X$, then it is of one of the following three
types:

(i) $E$ is decomposable and given by $L_a \oplus L_b$.

(ii) $E$ is not decomposable and an extension of line bundle; in this case, we can write
$E$ as $L_a \otimes E'$, where $E'$ is a non-trivial extension of the form:
\[ 0 \rightarrow L_{\mu_1^{m_1}\dots\mu_n^{m_n}} \rightarrow E' 
\rightarrow \mathcal{O} \rightarrow 0, \]
with $m_1,\dots,m_n$ non-negative integers.

(iii) $E$ is not an extension of line bundle; 
in this case, we can write $E$ as $L_a \otimes E'$, 
where $E'$ is the {\em unique} locally-free extension of the form:
\[ 0 \rightarrow  L_{\mu_1^{m_1} \dots \mu_i^{-k_i} \dots \mu_j^{-k_j} \dots \mu_n^{m_n}} 
\rightarrow E' \rightarrow I_{k_i k_j} \rightarrow 0, \]
with non-negative integers $m_1, \dots, m_n$ that are not all zero.
\end{proposition}

\begin{remark}\label{length-2 rank-2 projective resolutions}
Note that if $n \geq 4$, then the $H_{k_i k_j}$'s are the only codimension 2
subvarieties $Z$ of $X$ that are locally complete intersections   
whose ideal sheaves $I_Z$ admit projective resolutions
of the form $0 \rightarrow L \rightarrow V \rightarrow I_Z \rightarrow 0$,
where $L$ is a line bundle and $V$ is a rank-2 vector bundle on $X$. 
When $n=3$, one also has 
$Z = H_{k_i k_j} + H_{k_i k_l}$, with $k \neq l$.
\end{remark}

\begin{proof}
The proof of (i) and (ii) is the same as that of Proposition
\ref{classification - trivial Chern class}.
We therefore assume that $E$ is not an extension of line bundles. 
Since all vector bundles on $X$ are filtrable, it can then be written
in the form $E = L_a \otimes E'$, where $E'$ is a locally free extension of the form
\begin{equation}\label{general extension} 
0 \rightarrow \mathcal{O} \rightarrow E' \rightarrow L_\beta \otimes I_Z \rightarrow 0,
\end{equation}
$I_Z$ is the ideal sheaf of a codimension 2 subvariety $Z$ of $X$
that is a locally complete intersection, and $a,\beta \in \mathbb{C}^\ast$. 
By Remark \ref{length-2 rank-2 projective resolutions}, 
we must have $Z = H_{k_ik_j}$ or $H_{k_ik_{j_1}} + H_{k_ik_{j_2}}$,
for some $k_i,k_{j_1},k_{j_2} \geq 1$.
We therefore have to determine which extensions of the form \eqref{general extension}
give rise to non-isomorphic rank-2 vector bundles:
this is done by analysing the exact sequence
\begin{equation}\label{exact sequence}
0 \rightarrow H^1(X,L_{\beta^{-1}}) \rightarrow {\rm Ext}^1(X;I_Z,L_{\beta^{-1}})
\rightarrow H^0(X,Ext^1(I_Z,L_{\beta^{-1}})) \rightarrow H^2(X,L_{\beta^{-1}}). 
\end{equation} 

Recall that an extension class $\xi \in {\rm Ext}^1(X;I_Z,L_{\beta^{-1}})$ determines
a locally free sheaf if and only if its image in $H^0(X,Ext^1(I_Z,L_{\beta^{-1}}))$
generates the stalk of $Ext^1(I_Z,L_{\beta^{-1}})$ at every point of $X$.
Since $Ext^1(I_Z,L_{\beta^{-1}})$ is a torsion sheaf supported on 
each $H := H_{k_ik_j}$ appearing in $Z$, 
a necessary condition for the existence of a
locally free extension is therefore that 
$h^0(H,Ext^1(I_Z,L_{\beta^{-1}})|_H) \neq 0$ for each $H$. 
Note that the restriction of $Ext^1(I_Z,L_{\beta^{-1}})$ to $H$ can be identified with 
the line bundle $\det(\mathcal{N}_{H/X}) \otimes L_{\beta^{-1}}$, 
which is isomorphic to $L_{\beta^{-1}\mu_i^{k_i}\mu_j^{k_j}}$ 
(because $\mathcal{N}_{H/X}^\ast \cong I_H/I_H^2$
is generated as an $\mathcal{O}_H$-module by $z_i^{k_i}$ and $z_j^{k_j}$).
\vspace{.05in}

Let us first assume that $n=3$.
Note that in this case the underlying space of $H$ is 
the elliptic curve $\mathbb{C}^\ast / \mu_k$ with $k \neq i,j$;
consequently, a global section of $L_{\beta^{-1}\mu_i^{k_i}\mu_j^{k_j}}$
on $H$ is given by a holomorphic function $g(z)$ such that 
\[ g(\mu z) = \beta^{-1}\mu_i^{k_i}\mu_j^{k_j}g(z) \]
whose Laurent series
expansion is of the form
$\sum_{t = -\infty}^\infty \sum_{s = 0}^{k_j - 1} \sum_{r = 0}^{k_i - 1}
a_{rst}z_i^rz_j^sz_k^t$. 
One can easily verify that such a function exists if and only if
\[ \beta = \mu_i^{m_i}\mu_j^{m_j}\mu_k^{-\nu}, \]
where $m_i,m_j,\nu$
are integers such that $1 \leq m_l \leq k_l$ for $l = i,j$,
so that 
\[ g(z) = a_0 z_i^{k_i - m_i}z_j^{k_j - m_j}z_k^{-\nu}\] 
for some $a_0 \in \mathbb{C}$.
Thus, $h^0(H,Ext^1(I_Z,L_{\beta^{-1}})|_H) \neq 0$ if and only if
$\beta$ satisfies \eqref{condition on beta}, in which case
$h^0(H,Ext^1(I_Z,L_{\beta^{-1}})|_H) = 1$.
However, the function $g(z)$ must also generate the stalk of 
$Ext^1(I_Z,L_{\beta^{-1}})$ at every point,
implying that we must have $m_i = k_i$ and $m_j = k_j$.
A necessary condition for the existence of a locally free extension is therefore that
\begin{equation}\label{condition on beta} 
\beta = \mu_i^{k_i}\mu_j^{k_j}\mu_k^{-\nu},
\end{equation}
for some integer $\nu$.

Let us first consider the case $Z = H_{k_ik_j}$. Referring to 
\eqref{condition on beta}, this means that $\beta = \mu_i^{k_i}\mu_j^{k_j}\mu_k^{-\nu}$.  
We have to determine which values of $\nu$
give rise non-isomorphic rank-2 vector bundles.
We shall see that there are two different cases depending on whether or not 
$\nu$ is positive.
If $\nu > 0$, then $h^q(X,L_{\beta^{-1}}) = 0$ for $q=1,2$, so that
\eqref{exact sequence} reduces to:
\[ {\rm Ext}^1(I_Z,L_{\beta^{-1}}) = H^0(X,Ext^1(I_Z,L_{\beta^{-1}})) 
= \mathbb{C}.\]
Therefore, the extension determined by a non-zero element of 
${\rm Ext}^1(I_Z,L_{\beta^{-1}})$ determine,
up to isomorphism, a unique rank-2 vector bundle.
Note that this bundle cannot be an extension of line bundles.

If $\nu \leq 0$, although $h^1(X,L_{\beta^{-1}}) = 0$, we have
$h^2(X,L_{\beta^{-1}}) \neq 0$.
Let us assume that a class $\xi$ in 
${\rm Ext}^1(I_Z,L_{\beta^{-1}})$ generates a rank-2 vector bundle $E'$. 
In this case, $\mathcal{O}$ is not a destabilising line bundle of $E'$. 
Indeed, given the exact sequence
\begin{equation}\label{projective resolution} 
0 \rightarrow L_{\mu_i^{-k_i}\mu_j^{-k_j}} \rightarrow 
L_{\mu_i^{-k_i}} \oplus L_{\mu_j^{-k_j}} \rightarrow I_H \rightarrow 0, 
\end{equation}
one verifies that $L_{\mu_i^{k_i}\mu_k^{-\nu}}$ and 
$L_{\mu_j^{k_j}\mu_k^{-\nu}}$
not only map non-trivially to $L_{\mu_i^{k_i}\mu_j^{k_j}\mu_k^{-\nu}} \otimes I_H$, 
but also to $E'$; 
furthermore, they are the line bundles of maximal degree mapping into $E'$. 
Consequently, the quotient sheaves 
$E' / L_{\mu_i^{k_i}\mu_k^{-\nu}}$ and 
$E' / L_{\mu_j^{k_j}\mu_k^{-\nu}}$ are torsion free, 
so that, for example, we have
an exact sequence
\[0 \rightarrow \mathcal{O} \rightarrow L_{\mu_i^{-k_i}\mu_k^\nu} \otimes E'
\rightarrow L_{\mu_i^{-k_i}\mu_j^{k_j}\mu_k^\nu} \otimes I_{Z'} \rightarrow 0, \]
where $Z'$ is either empty or a locally complete intersection of codimension 2 in $X$.
If $\nu=0$, then $E' = L_{\mu_i^{k_i}} \oplus L_{\mu_j^{k_j}}$,
implying that $Z'$ is empty.
However,  if $\nu<0$, then $\beta' = \mu_i^{-k_i}\mu_j^{k_j}\mu_k^\nu$
does not satisfy condition \eqref{condition on beta}, 
leading to a contradiction. 
Hence, we get locally free sheaves only for $\nu > 0$, 
in which case they are not extensions of line bundles.

Let us now assume that $Z = H_{k_ik_{j_1}} + H_{k_ik_{j_1}}$.
In this case, since $\beta$ must satisfy \eqref{condition on beta} for both
$(k_i,k_{j_1})$ and $(k_i,k_{j_2})$, we have 
$\beta = \mu_i^{k_i}\mu_j^{k_{j_1}}\mu_k^{k_{j_2}}$.
As above, one can show that $E'$ decomposes as 
$L_{\mu_i^{k_i}} \oplus L_{\mu_j^{k_{j_1}}\mu_k^{k_{j_2}}}$.
The only extensions giving rank-2 vector bundles that are not extensions of line 
bundles therefore correspond to codimension 2 subvarieties of the form $Z = H_{k_ik_j}$. 
This proves (iii) for $n=3$.

Finally, let us assume that $n\geq 4$ so that the underlying space of
$H$ is a generic Hopf manifold of dimension at least 2. 
As in the case $n=3$, we can easily show that
there exists a locally free extension $E'$ if and only if
\[\beta = \mu_1^{-m_1} \dots \mu_i^{k_i} \dots \mu_j^{k_j} \dots \mu_n^{-m_n},\] 
for non-negative integers $m_1, \dots, m_n$;
note that the $m_l$'s must all be non-negative because $H$ 
is now a Hopf manifold. 
Moreover, $E'$ is not an extension of line bundles,  
unless all the $m_l$'s are zero, in which case 
$E' = L_{\mu_i^{k_i}} \oplus L_{\mu_j^{k_j}}$. 
\end{proof}

We finish the paper by determining stability conditions for
holomorphic rank-2 vector bundles on higher dimensional Hopf manifolds. 

\begin{theorem}\label{stable bundles - generic manifolds}
Let $E$ be a rank-2 vector bundle on a generic Hopf manifold $X$ of dimension
greater than 2.

(i) If $E$ is an extension of line bundles, then it is unstable. 

(ii) Otherwise, $E = L_a \otimes E'$
for $a \in \mathbb{C}^\ast$ 
and a non-trivial extension $E'$ of 
$I_{H_{{k_i k_j}}}$ by $L_{\mu_1^{m_1} \dots \mu_i^{-k_i} \dots \mu_j^{-k_j} \dots \mu_n^{m_n}}$,
where the non-negative integers $m_1, \dots, m_n$ are not all zero 
(see Proposition \ref{classification bundles generic manifold} (iii)).
In this case, $E$ is stable if and only if 
\[ \prod_{\stackrel{0 \leq l \leq n}{l \neq i,j}} | \mu_l^{m_l} | 
> | \mu_i^{k_i} \mu_j^{k_j}|. \] 
Note that one can always find integers $m_1, \dots, m_n$ that satisfy this equation.
Consequently, stable rank-2 vector bundles exist on $X$.
\hfill \qedsymbol
\end{theorem}
\begin{proof}
We begin by determining the destabilising line bundles of rank-2
vector bundles; given that decomposable bundles are automatically unstable,
we only consider the indecomposable case.
Let $E$ be an indecomposable rank-2 vector bundle on $X$.
Suppose that $E = L_a \otimes E'$, 
where $a \in \mathbb{C}^\ast$ and $E'$ is one of the non-trivial extension 
$0 \rightarrow L \rightarrow E' \rightarrow I_Z \rightarrow 0$
of Proposition \ref{classification bundles generic manifold} (ii) or (iii).
Therefore, any line bundle mapping non-trivially into $E$ is of the
form 
\[ L \otimes L_{a \mu_1^{-l_1} \dots \mu_n^{-l_n}},\] 
for non-negative integers $l_1, \dots, l_n$.
This is obvious when $Z$ is empty.
Otherwise, if $Z = H_{{k_i k_j}}$, then this comes
from the fact that every bundle mapping to $I_{H_{{k_i k_j}}}$ also
maps to $L$.
Consequently, since $L \otimes L_a$ is the unique destabilising line bundle of $E$,
the theorem follows from the definition of degree. 
\end{proof}


\begin{thebibliography}{References.}
\bibitem[BPV]{BPV} W. Barth, C. Peters, and A. Van de Ven, 
{\em Compact complex surfaces}, Springer-Verlag
Berlin, Heidelberg, New York, 1984. 

\bibitem[Bo]{Bottacin-Poisson} F. Bottacin, 
{\em Poisson structures on moduli spaces of 
sheaves over Poisson surfaces}, Invent. Math.  {\bf 121} (1995) 421-436.

\bibitem[B1]{Braam1} P. J. Braam,
{\em Magnetic monopoles on 3-manifolds}, J. Diff. Geom. {\bf 30} 425-464.

\bibitem[B2]{Braam2} P. J. Braam,
{\em A Kaluza-Klein approach to hyperbolic 3-manifolds}, Enseignement Math. (2) {\bf 34} (1989)
275-312.

\bibitem[BH]{B-H} P. J. Braam and J. Hurtubise,
{\em Instantons on Hopf surfaces and monopoles on solid tori}, J. Reine Angew. Math. {\bf 400} (1989)
146-172.

\bibitem[Br]{Brinzanescu} V. Br\^{\i}nz\u{a}nescu,
{\em Double covers and vector bundles}, 
An.\ Stiint.\ Univ.\ Ovidius\ Constanta,\ Ser.\ Mat.\ {\bf 9} (2001) no 1 21-26. 

\bibitem[BrMo1]{Brinzanescu-Moraru} V. Br\^{\i}nz\u{a}nescu and R. Moraru,
{\em Holomorphic rank-2 vector bundles on non-K\"{a}hler elliptic surfaces}, preprint
arXiv:math.AG/0306191.

\bibitem[BrMo2]{Brinzanescu-Moraru1} V. Br\^{\i}nz\u{a}nescu and R. Moraru,
{\em Twisted Fourier-Mukai transforms and bundles on non-K\"ahler elliptic surfaces}, 
preprint arXiv:math.AG/0309031.

\bibitem[BrMo3]{Brinzanescu-Moraru2} V. Br\^{\i}nz\u{a}nescu and R. Moraru,
{\em Stable bundles on non-K\"ahler elliptic surfaces}, 
to appear in Comm.\  Math.\  Phys.;
preprint arXiv:math.AG/0306192.

\bibitem[Bh]{Buchdahl} N. P. Buchdahl, 
{\em Hermitian-Einstein connections and stable
vector bundles over compact complex surfaces}, Math. Ann. {\bf 280} (1988) 625-648.

\bibitem[D]{D}
R. Donagi, 
{\em Principal bundles on elliptic fibrations},
Asian\ J.\ Math.\ {\bf 1(2)} (1997) 214-223.

\bibitem[F]{F} 
R. Friedman,
{\em Algebraic Surfaces and Holomorphic Vector
Bundles}, UTX, Springer, New York Berlin Heidelberg, 1998.

\bibitem[FM]{FM}
R. Friedman \and J. W. Morgan, 
{\em Smooth Four-Manifolds and Complex Surfaces}, Springer-Verlag, 1994.

\bibitem[FMW]{FMW} 
R. Friedman, J. Morgan, and E. Witten, 
{\em Vector bundles over
elliptic fibrations}, J. Algebraic Geom. {\bf 2} (1999) 279-401.

\bibitem[G]{Gauduchon} P. Gauduchon, 
{\em Le th\'{e}or\`{e}me de l'excentricit\'{e} nulle},
C. R. A. S. Paris {\bf 285} (1977) 387-390.

\bibitem[LT]{LT} M. L\"ubke and A. Teleman, 
{\em The Kobayashi-Hitchin correspondence}, World
   Scientific Publishing Co., Inc., River Edge, NJ, 1995.

\bibitem[Ma1]{Mall-line bundles}
D. Mall,
{\em The cohomology of line bundles on Hopf manifolds},
Osaka J. Math. {\bf 28} (1991) 999-1015.

\bibitem[Ma2]{Mall-filtrable}
D. Mall,
{\em On holomorphic vector bundles on Hopf manifolds with trivial pullback
on the universal covering},
Math. Ann. {\bf 294} (1992) 719-740.

\bibitem[Ma3]{Mall-rank-2}
D. Mall,
{\em Filtrable bundles of rank 2 on Hopf surfaces},
Math. Z. {\bf 212} (1993) 239-256. 

\bibitem[Mo1]{Moraru} 
R. Moraru, {\em Integrable systems
associated to a Hopf surface}, Canad. J. Math. {\bf 55} (3) (2003) 609-635.



\bibitem[T]{T} A. Teleman, {\em Moduli spaces of stable bundles on non-K\" ahler elliptic
fibre bundles over curves}, Expo. Math. \textbf{16} (1998) 193-248.

\bibitem[V1]{V1} M. Verbitsky, 
{\em Stable bundles on positive principal elliptic fibrations}, preprint arXiv:math.AG/0403430.

\bibitem[V2]{V2} M. Verbitsky, 
{\em Holomorphic bundles on diagonal Hopf manifolds}, math.AG/0408391.
\end{thebibliography}
\end{document}